\documentclass[a4paper, 12pt, oneside, notitlepage]{amsart}
\usepackage[margin=3cm]{geometry}
\usepackage{amsmath,amssymb,amsthm,graphicx,mathrsfs,bbm,url}
\usepackage{amsthm}
\usepackage{wrapfig}
\usepackage{enumitem}
\usepackage{mathtools}
\usepackage[utf8]{inputenc}
 \usepackage[T1]{fontenc}
\usepackage[usenames,dvipsnames]{color}
\usepackage[colorlinks=true,linkcolor=Red,citecolor=Green]{hyperref}
\usepackage[super]{nth}
\usepackage[open, openlevel=2, depth=3, atend]{bookmark}
\hypersetup{pdfstartview=XYZ}
\usepackage[font=footnotesize]{caption}
\usepackage{a4wide}

\usepackage{epstopdf}
 
\usepackage{hyperref}

\theoremstyle{plain}
\newtheorem{theorem}{Theorem}[section]
\newtheorem*{theorem*}{Theorem}

\newtheorem{proposition}[theorem]{Proposition}

\theoremstyle{definition}

\theoremstyle{remark}

\numberwithin{equation}{section}

\newcommand{\R}{\mathbb{R}}

\newcommand{\Z}{\mathbb{Z}}

\newcommand{\T}{\mathbb{T}}

\newcommand{\HH}{\mathbb{H}}

\newcommand{\eps}{\varepsilon}

\newcommand{\mc}{\mathcal}

\DeclareMathOperator{\vol}{vol}

\DeclareMathOperator{\Op}{Op}

\DeclareMathOperator{\Leb}{Leb}

\newcommand{\be}{\begin{equation}}
\newcommand{\ee}{\end{equation}}

% ...

% ...

\title
[The magnetic Laplacian on hyperbolic surfaces]
{Semiclassical analysis of the magnetic Laplacian on hyperbolic surfaces}

\author{Thibault Lefeuvre}

\address{Université de Paris and Sorbonne Université, CNRS, IMJ-PRG, F-75006 Paris, France.}

\email{tlefeuvre@imj-prg.fr}

\begin{document}

\begin{abstract}
The magnetic Laplacian on hyperbolic surfaces provides a rich analytic framework in
which a variety of quantum phenomena emerge. The present note, written for the \emph{Proceedings of the
Journées EDP 2025}, is a concise overview of the main results obtained in \cite{Charles-Lefeuvre-25,Chabert-Lefeuvre-26,Charles-Lefeuvre-26}. 
\end{abstract}

\maketitle

\section{Introduction}

\subsection{Overview of the results}

We are concerned with the semiclassical behaviour of eigenfunctions
\[
k^{-2}\Delta_k u_k=(E+\varepsilon_k)u_k,
\qquad
u_k\in C^\infty(\Sigma,L^{\otimes k}),
\]
where $L\to\Sigma$ is a Hermitian line bundle over a closed oriented hyperbolic surface
$\Sigma$ of genus $\geq2$, equipped with a unitary connection $\nabla$, and $k \geq 0$ is an integer. The induced connection on $L^{\otimes k}$ has curvature
$F_{\nabla^k}=-ikB\,\mathrm{vol}$, where $\vol$ is the Riemannian volume form, and throughout the paper the magnetic field
$B$ is assumed to be \emph{constant}. The results discussed in this note can be summarized as follows:

\begin{itemize}
\item When the energy lies below the critical value $E<E_c$, the associated classical
Hamiltonian flow is completely periodic. In this regime, one can show that every
flow-invariant probability measure arises as a semiclassical defect measure (Theorem \ref{theorem:main}, item (i)).
Moreover, the Hörmander bound for the $L^\infty$ norm of eigenfunctions is saturated
by explicit families of states, which we call \emph{magnetic zonal states} and may be viewed as analogues of zonal harmonics
on the sphere. They exhibit concentration on a two-dimensional torus in phase space (Theorem \ref{theorem:Linfty}, item (i)). However, unlike zonal harmonics, their $L^\infty$ norm only blows up at a single point due to the absence of ``magnetic'' conjugate points.

\item At the critical energy $E=E_c$, the picture changes abruptly. The magnetic
Hamiltonian flow becomes conjugate to the horocyclic flow, which is uniquely
ergodic. As a consequence, the Liouville measure is the only possible semiclassical
defect measure (Theorem \ref{theorem:main}, item (ii)), and a form of quantitative unique ergodicity holds. More precisely,
one obtains a polynomial convergence rate of expectation values
$\langle \Op_k(a)u_k,u_k\rangle_{L^2}$ towards their Liouville averages (Theorem \ref{theorem:horocyclic}). This
quantitative control can in turn be exploited to derive a polynomial improvement
over Hörmander’s $L^\infty$ bound for eigenfunctions at the critical energy (Theorem \ref{theorem:Linfty}, item (ii)).

\item Above the critical threshold $E>E_c$, the dynamics enters a genuinely hyperbolic
regime: the magnetic flow is conjugate to a reparametrization of the geodesic flow on
the unit tangent bundle, which is Anosov. The analysis of eigenfunctions in this case
falls within the scope of the Quantum Unique Ergodicity conjecture for
negatively-curved Riemannian manifolds \cite{Anantharaman-Nonnenmacher-07,Anantharaman-08,Dyatlov-Jin-18,Dyatlov-Jin-Nonnenmacher-22}, a central and largely open problem in quantum
chaos. Namely, there is a subset of eigenfunctions of density $1$ which equidistribute towards the Liouville measure (Theorem \ref{theorem:main}, item (iii)) and it is conjectured that no exceptional subsequences should exist.
\end{itemize}

\subsection{Setup}

Let $(\Sigma,g)$ be a compact, connected, oriented surface of genus $\textsl{g}\geq2$,
equipped with a hyperbolic metric of constant curvature $-1$.
Let $L\to\Sigma$ be a Hermitian complex line bundle endowed with a unitary connection
$\nabla$. Its curvature is given by
\[
F_\nabla=-i\,B\,\mathrm{vol}\in C^\infty(\Sigma,\Lambda^2T^*\Sigma),
\]
where $B\in C^\infty(\Sigma)$ and $\mathrm{vol}$ denotes the Riemannian volume form.
The function $B$ will be referred to as the \emph{magnetic field}.
The associated \emph{magnetic Laplacian} is defined by
\begin{equation}
\Delta_L:=\tfrac12\nabla^*\nabla:
C^\infty(\Sigma,L)\longrightarrow C^\infty(\Sigma,L).
\end{equation}
For any integer $k\geq0$, we consider the tensor power $L^{\otimes k} \to \Sigma$ endowed with the
induced connection $\nabla^{\otimes k}$, and define in the same manner a Laplacian
$\Delta_k$ acting on sections of $L^{\otimes k}$. After rescaling by $k^{-2}$,
the operator $k^{-2}\Delta_k$ enters the semiclassical regime with parameter
$h:=k^{-1}$.

Throughout this work, we further assume that the magnetic field $B$ is constant.
By the Gauss--Bonnet theorem, this implies the quantization condition
$2B(\textsl{g}-1)\in\mathbb{Z}$.
Our goal is to study the semiclassical behaviour of normalized eigenstates
\begin{equation}
\label{equation:eigenstates}
k^{-2}\Delta_k u_k=(E+\varepsilon_k)u_k,
\qquad
u_k\in C^\infty(\Sigma,L^{\otimes k}),
\quad
\|u_k\|_{L^2}=1,
\end{equation}
in the limit $k\to+\infty$, where $E\geq0$ is fixed and $\varepsilon_k\to0$.

We say that a sequence $(u_k)_{k \geq 0}$ satisfying \eqref{equation:eigenstates} converges to a \emph{semiclassical defect measure}
$\mu$ on $T^*\Sigma$ if for every symbol $a\in C^\infty(T^*\Sigma)$,
\begin{equation}
\label{equation:convergence}
\langle \Op_k(a)u_k,u_k\rangle_{L^2}
\longrightarrow_{k \to +\infty}
\int_{T^*\Sigma}a(x,\xi)\,d\mu(x,\xi),
\end{equation}
where $\Op_k$ denotes a magnetic semiclassical quantization.
We write $u_k\rightharpoonup\mu$ for this convergence.
The normalization of $u_k$ implies that $\mu$ is a probability measure, and
its support is contained in the compact energy hypersurface $\{p=E\}\subset T^*\Sigma$,
where
\[
p(x,\xi):=\tfrac12|\xi|_g^2
\]
is the principal symbol of $k^{-2}\Delta_k$.
Although the limit \eqref{equation:convergence} may fail to exist along the full
sequence, it always exists along subsequences, which we do not always relabel.
As we shall see, the structure of possible defect measures depends crucially on
the value of $E$, leading to three distinct semiclassical regimes.

\subsection{Magnetic flow in phase space} Let $\omega_0$ denote the canonical Liouville symplectic form on $T^*\Sigma$, and set
\begin{equation}
\label{equation:twisted-form}
\Omega:=\omega_0+i\,\pi^*F_\nabla,
\end{equation}
where $\pi:T^*\Sigma\to\Sigma$ is the canonical projection.
The form $\Omega$ is symplectic and incorporates the magnetic contribution.
The Hamiltonian vector field $H_p^\Omega$ associated with $p$ is defined by
\[
dp(\bullet)=\Omega(\bullet,H_p^\Omega),
\]
and the corresponding Hamiltonian flow $(\Phi_t)_{t\in\mathbb{R}}$ is called the
\emph{magnetic flow}.
Each energy hypersurface $\{p=E\}$ is invariant under this flow and carries a natural
smooth invariant probability measure, the \emph{Liouville measure},
denoted $\mu_{\mathrm{Liouv}}$.

Let $S\Sigma$ be the unit tangent bundle of $\Sigma$.
We denote by $(\varphi_t)_{t\in\mathbb{R}}$ the geodesic flow on $S\Sigma$, by $(R_t)_{t\in\mathbb{R}}$ the $2\pi$-periodic rotation
in the fibers of $S\Sigma$, and by $(h_t)_{t\in\mathbb{R}}$ the stable horocyclic flow. Note that $S\Sigma = S\HH^2/\Gamma = \mathrm{PSL}(2,\R)/\Gamma$ for some lattice $\Gamma < \mathrm{PSL}(2,\R)$. The three flows defined above are respectively generated by the following elements in $\mathfrak{sl}(2,\R)$, seen as right-invariant vector fields on $\mathrm{PSL}(2,\R)$:
\[
X = \begin{pmatrix} 1/2 & 0 \\ 0 & -1/2 \end{pmatrix}, \qquad V =\begin{pmatrix} 0 & 1/2 \\ -1/2 & 0\end{pmatrix}, \qquad U_+ =  \begin{pmatrix} 0 & 1 \\ 0 & 0\end{pmatrix}.
\]
Define the \emph{critical energy}
\[
E_c:=\tfrac12 B^2,
\]
and set
\[
T_E=
\begin{cases}
(B^2-2E)^{-1/2}, & E<E_c,\\[0.3em]
(2E-B^2)^{-1/2}, & E>E_c.
\end{cases}
\]
It is well known (see \cite[Proposition 2.5]{Charles-Lefeuvre-25} for instance) that for $E>0$,
the magnetic flow on $\{p=E\}$ is smoothly conjugate to:
\begin{enumerate}[label=(\roman*)]
\item the reparametrized rotation flow $(R_{t/T_E})_{t \in \R}$ when $E<E_c$ (elliptic case);
\item the horocyclic flow $(h_t)_{t \in \R}$ when $E=E_c$ (parabolic case);
\item the reparametrized geodesic flow $(\varphi_{t/T_E})_{t \in \R}$ when $E>E_c$ (hyperbolic case).
\end{enumerate}
This is a straightforward consequence of the fact that the magnetic flow on the energy shell $\{p=E\}$ is conjugate to the flow generated by $\sqrt{2E} X - BV$ on $S\Sigma$. Then, given $Y \in \mathfrak{sl}(2,\R)$, it is conjugate to one of the above flows depending on the trichotomy $\det Y <  0, \det Y > 0$ or $\det Y = 0$.

\section{Quantum limits}

\subsection{Semiclassical defect measures}

The dynamical transition at the classical level for the magnetic flow has a direct quantum counterpart,
summarized in the following theorem.

\begin{theorem}[Three semiclassical regimes]
\label{theorem:main}
Under the above assumptions on $(\Sigma,g)$ and $(L,\nabla)$, the following statements hold:
\begin{enumerate}[label=\emph{(\roman*)}]
\item \textbf{\emph{Low-energy regime.}}
If $0\leq E<E_c$, then for every probability measure $\mu$ on $\{p=E\}$ that is invariant
under the magnetic flow $(\Phi_t)_{t\in\mathbb{R}}$, there exists a sequence
$(u_k)_{k\geq0}$ satisfying~\eqref{equation:eigenstates} such that
\[
u_k\rightharpoonup_{k\to\infty}\mu .
\]

\item \textbf{\emph{Critical-energy regime.}}
If $E=E_c$, then for any sequence $(u_k)_{k\geq0}$ satisfying
\eqref{equation:eigenstates}, the associated semiclassical defect measure is the Liouville measure:
\[
u_k\rightharpoonup_{k\to\infty}\mu_{\mathrm{Liouv}} .
\]

\item \textbf{\emph{High-energy regime.}}
Fix $E_c<a<b$. For each $k\geq0$, let $(u_{k,j})_{j\in J_k}$ be an orthonormal family of
eigenfunctions of $k^{-2}\Delta_k$ with eigenvalues $\lambda_{k,j}\in[a,b]$. Then there
exists a subset $J_k^\ast\subset J_k$ of density one such that, for any sequence
$(u_{k_n,j_n})$ with $j_n\in J_{k_n}^\ast$ and $k_n\to+\infty$, one has
\[
u_{k_n,j_n}\rightharpoonup_{n\to\infty}\mu_{\mathrm{Liouv}} .
\]
\end{enumerate}
\end{theorem}

In the high-energy regime $E>E_c$, the spectrum of the
magnetic Laplacian is directly related to that of the Laplace--Beltrami operator on
$(\Sigma,g)$. In this setting, it is
conjectured that the Liouville measure is the only possible semiclassical defect
measure, in accordance with the Quantum Unique Ergodicity (QUE) conjecture, see \cite{Dyatlov-22} and the references therein for further discussion.

The low-energy regime $E<E_c$ is more subtle. It relies on the construction of
eigenstates that concentrate microlocally along prescribed periodic orbits of the
magnetic flow. This construction uses an averaging argument in the spirit of
Weinstein~\cite{Weinstein-77}, see \S\ref{ssection:weinstein}.

The critical regime $E=E_c$ provides a particularly transparent, albeit somewhat
degenerate, manifestation of QUE. Indeed, any semiclassical defect measure must be
invariant under the magnetic Hamiltonian flow, which coincides on $\{p=E_c\}$ with the
horocyclic flow. Since the horocyclic flow is uniquely ergodic
\cite{Furstenberg-73}, the Liouville measure is the only invariant probability measure,
and hence the only possible defect measure.

In the low-energy regime, we also obtain more general results regarding the perturbed operator $k^{-2}(\Delta_k + V)$ where $V \in C^\infty(\Sigma)$ is a smooth potential; see \cite[Theorem 1.3]{Charles-Lefeuvre-25} for a detailed discussion. From a geometric viewpoint, perturbing $\Delta_k$ by a potential is, to some extent, analogous to perturbing the round metric on the sphere to a Zoll metric. In the forthcoming article \cite{Charles-Lefeuvre-26}, the operator $k^{-2}(\Delta_k + V)$ is studied on surfaces with negative sectional curvature (not necessarily hyperbolic) in the regime where $E \simeq 0$.

\subsection{Quantitative Quantum Unique Ergodicity}

In the critical
regime, we can in fact obtain a quantitative refinement of this convergence. Let $0<\theta<\tfrac12$ be such that
$\theta(1-\theta)\leq\lambda_1(\Sigma)$, where $\lambda_1(\Sigma)$ denotes the first
non-zero eigenvalue of the Laplace--Beltrami operator on functions. Denote by
$S_c^\ast\Sigma := \{p=E_c\}$ the critical energy hypersurface. 

\begin{theorem}[Polynomial convergence rate at the critical energy]
\label{theorem:horocyclic}
Assume that~\eqref{equation:eigenstates} holds with $E=E_c$ and that
$\varepsilon_k\leq h^\ell$ for some $0 < \ell \leq 1/15$. Then there exists a constant
$C=C(\ell)>0$ such that, for any symbol $a\in C^\infty(T^\ast\Sigma)$ supported in
$\{p\leq10E_c\}$, one has
\[
\left|
\langle \Op_k(a)u_k,u_k\rangle_{L^2}
-
\int_{S_c^\ast\Sigma} a\,d\mu_{\mathrm{Liouv}}
\right|
\leq
C\,k^{-\theta \ell/4100}\,
\|a\|_{C^{17}(T^\ast\Sigma)} .
\]
\end{theorem}

The exponent appearing in the remainder term is not expected to be optimal, and no
attempt has been made to improve it. Of course, one can take $\ell \geq 1/15$, but then the exponent of convergence is \emph{a priori} not better than $\theta \ell (15\cdot 4100)^{-1}$.

The compact support assumption on $a$ may be
relaxed: one may equivalently consider symbols $a\in S^m(T^\ast\Sigma)$, using the fact
that the eigenstates $u_k$ concentrate on the energy hypersurface $\{p=E_c\}$. This
modification only produces a negligible $\mathcal{O}(k^{-\infty})$ error.

Finally, we mention that related quantitative results in the direction of QUE were
recently obtained by Morin and Rivière~\cite{Morin-Riviere-25} in a different
semiclassical regime, where the magnetic field is constant and the semiclassical
parameter is the inverse square root of the energy, in the setting of magnetic
Laplacians on the torus.

\subsection{Proof ideas} \label{ssection:weinstein}

The first eigenvalues of $\Delta_k = \tfrac{1}{2}\nabla_k^*\nabla_k$ acting on $C^\infty(\Sigma,L^{\otimes k})$ are explicit and given for $0 \leq m < N_k := \lfloor kB\rfloor $ by
\begin{equation}
\label{equation:ev}
\lambda_{k,m} = kB (m+\tfrac{1}{2})-\dfrac{m(m+1)}{2}.
\end{equation}
As $k \rightarrow \infty$, $\lambda_{k, N_k - 1 } \sim k^2 E_c$ where $E_c
= \tfrac{1}{2} B^2$ is the critical energy. See \cite[Proposition 2.8]{Charles-Lefeuvre-25}.

\subsubsection{Weinstein's periodic operator} Let $\Pi_{k,m}$ be the $L^2$-orthogonal projector of $C^\infty(M, L^{\otimes k})$ onto
 $\ker(\Delta _k -\lambda_{k,m})$, and set:
\begin{equation}
\label{equation:a}
\mathbf{A}_k := k^{-1} \sum_{m =0 }^{ \lfloor kB \rfloor -1}  m \Pi_{k,m}.
\end{equation}
It can be established that $\mathbf{A}_k$ is a (twisted) compactly supported semiclassical operator with support in the unit disk bundle $D^*\Sigma := \{(x,\xi) \in T^*\Sigma ~|~ |\xi| < B\}$, see \cite[Lemma 4.1]{Charles-Lefeuvre-25}. 
Furthermore, by \eqref{equation:ev}, it can be verified that on the space
\[
\mc{I}_k := \bigoplus_{m=0}^{ \lfloor kB \rfloor -1}
\ker(\Delta_k-\lambda_{k,m}),
\]
the operator $\mathbf{A}_k$ satisfies the algebraic identity
\begin{equation}
\label{equation:relation}
k^{-2} \Delta_k = B(\mathbf{A}_k+ \tfrac{1}{2} k^{-1} )- \tfrac{1}{2} 
\mathbf{A}_k(\mathbf{A}_k+ k^{-1} ).
\end{equation}
Letting $a$ be the function defined on $D^*\Sigma$ through the relation
\[
p(x,\xi) = \beta (a ( x, \xi)), \qquad \forall (x,\xi) \in D^*\Sigma,
\]
where $\beta : [0,B] \rightarrow [0, \frac{1}{2} B^2 ]$ is the function $ \beta (s)  :=  B s  -\frac{1}{2} s^2$, a straightforward computation shows that $a$ is the principal symbol of $\mathbf{A}_k$.

The Hamiltonian flow generated by $a$ (computed with respect to the twisted symplectic $2$-form $\Omega$, see \eqref{equation:twisted-form}) is a $2\pi$-periodic reparametrization of the one generated by $p$. Observe that by construction $\mathrm{sp}(\mathbf{A}_k) \subset k^{-1} \Z_{\geq 0}$ and thus
\begin{equation}
\label{equation:cool}
e^{2i k \pi \mathbf{A}_k } = \mathbf{1}.
\end{equation}
Finally, note that
\begin{equation}
\label{eq:integrale_projecteur}
\Pi_{k,m} = \dfrac{1}{2\pi} \int_0^{2\pi} e^{-imt} e^{itk\mathbf{A}_k} d t.
\end{equation}
This turns out to be a Fourier Integral Operator in the semiclassical regime $k \to +\infty$ with an explicit canonical relation given by
\[
C = \{(\Phi_t(x,\xi) ; (x,\xi) ) ~:~ t \in [0,T_E], (x,\xi) \in \{p=E\}\}.
\]
See \cite{Weinstein-77,Zelditch-97} for related discussions.

\subsubsection{Gaussian beams} The spectral projector $\Pi_{k,m}$ can be used to produce magnetic Gaussian beams, namely sections that microlocally concentrate on a periodic bicharacteristic. We fix an energy $0 < E < E_c$. Given any $(x,\xi) \in \{p=E\} \subset T^*\Sigma$, its magnetic bicharacteristic $(\Phi_t(x,\xi))_{t \in \R}$ is $T_E$-periodic where $T_E = (B^2-2E)^{-1/2}$. One may associate a Gaussian wave packet
$\mathbf{e}_{k,x,\xi} \in C^\infty(\Sigma,L^{\otimes k})$ microlocally concentrated near
the phase-space point $(x,\xi)$.

For a given $k \geq 0$, we let $m_k$ be an integer such that $\lambda_{k,m_k} \to E$. By applying the spectral projector $\Pi_{k,m_k}$, one obtains a genuine eigenfunction
\[
\mathbf{f}_{k,x,\xi}
:= k^{1/4}\,\Pi_{k,m_k} \mathbf{e}_{k,x,\xi}
\in C^\infty(\Sigma,L^{\otimes k}),
\]
satisfying~\eqref{equation:eigenstates} at energy $E$.
This construction is detailed in \cite[Proposition~4.6]{Charles-Lefeuvre-25}. The resulting eigenstate
$\mathbf{f}_{k,x,\xi}$ is an $L^2$ (quasi-)normalized\footnote{There exists $C > 1$ such that $C^{-1} \leq \|\mathbf{f}_{k,x,\xi}\|_{L^2} \leq C$.} magnetic Gaussian beam, whose semiclassical
mass is concentrated along the periodic trajectory generated by $(x,\xi)$.

In particular, modulo renormalization, it converges to the semiclassical defect measure $\delta_\gamma$ which is the Dirac mass supported by the periodic bicharacteristic generated by $(x,\xi)$. As any measure on $\{p=E\}$ which is invariant by the magnetic flow $(\Phi_t)_{t \in \R}$ can be approximated by such Dirac masses, an elementary argument shows that any flow-invariant measure can be obtained as a semiclassical defect measure, thus proving Theorem \ref{theorem:main}, item (i).

\subsubsection{Quantitative QUE} 

The proof of Theorem \ref{theorem:horocyclic} relies on two key ingredients. The first is a quantitative version of the classical unique ergodicity for the horocyclic flow due to Burger \cite{Burger-90}, which asserts that for all $0 < \theta < 1/2$ such that $\theta(1-\theta) \leq \lambda_1(\Sigma)$ (the first non-zero eigenvalue of the Laplace-Beltrami operator on functions), the following estimate holds: there exists a constant $C > 0$ such that for all $T > 0$, for all $a \in H^3(S\Sigma)$,
\begin{equation}
\label{equation:burger}
\sup_{v \in S\Sigma} \left|\dfrac{1}{T} \int_0^T a(h_t(v)) d t - \int_{S\Sigma} a(v) d\mu_{\mathrm{Liouv}}(v)\right| \leq \dfrac{C\|a\|_{H^3(S\Sigma)}}{T^\theta}.
\end{equation}

The second ingredient is a long-time version of the Egorov theorem. In its standard form, the Egorov theorem asserts that for every observable $a \in C^\infty_{\mathrm{comp}}(T^*\Sigma)$:
\begin{equation}
\label{equation:egorov}
e^{itk^{-1}\Delta_k} \Op_k(a) e^{-itk^{-1}\Delta_k} = \Op_k(a \circ \Phi_t) + \mc{O}(k^{-1}),
\end{equation}
where the remainder is uniform with respect $t \leq T_k := \eps \log k$, provided $\eps > 0$ is small enough. The time $T_k$ is called the \emph{Ehrenfest time}. However, this is not enough to apply \eqref{equation:burger} as we want to go \emph{beyond} the Ehrenfest time and use it with polynomial times $t \simeq k^\delta$, where $\delta > 0$. Usually, this is only possible if the flow satisfies parabolic estimates of the form
\begin{equation}
\label{equation:parabolic}
\|f \circ \Phi_t\|_{C^n(T^*\Sigma)} \leq C_n \langle t \rangle^{m_n} \|f\|_{C^n(T^*\Sigma)}, \qquad \qquad \forall t \geq 0,
\end{equation}
for all functions $f \in C^\infty(T^*\Sigma)$, where $m_n \geq 0$. However, the estimate \eqref{equation:parabolic} fails in a drastic way on $T^*\Sigma$ as the magnetic flow $(\Phi_t)_{t \in \R}$ is hyperbolic on $\{p > E_c\}$, which implies that for all compact $K \subset T^*\Sigma$, $n \geq 0$, and $f \in C^\infty_{\mathrm{comp}}(T^*\Sigma)$ with support in $K$,
\[
\|f \circ \Phi_t\|_{C^n(T^*\Sigma)} \leq C_n e^{\gamma n t} \|f\|_{C^n(T^*\Sigma)}, \qquad \forall t \geq 0,
\]
and this estimate cannot be improved. The exponent $\gamma > 0$ depends on $K$; it is called a \emph{Lyapunov exponent} and is the maximal expansion rate of the flow.

The crucial observation is that, as $E \to E_c$, the hyperbolicity of the flow weakens in such a way that, in restriction to the critical energy shell $\{p=E_c\}$, the parabolic estimate \eqref{equation:parabolic} holds. More precisely, for all $n \geq 0$, there exists a constant $C_n > 0$ such that for all $f \in C^\infty_{\mathrm{comp}}(T^*\Sigma)$ with support in $\{p \leq E\}$:
\begin{equation}
\label{equation:good}
\|f \circ \Phi_t\|_{C^n(T^*\Sigma)} \leq C_n \langle t \rangle^{m_n} e^{(2(E-E_c)_+)^{1/2}nt} \|f\|_{C^n(T^*\Sigma)}, \qquad \forall t \geq 0,
\end{equation}
where $x_+ = \max(x,0)$.

As the eigenstates $u_k \in C^\infty(\Sigma,L^{\otimes k})$ satisfying \eqref{equation:eigenstates} with $E=E_c$ concentrate on the energy shell $\{p=E_c\}$, one can replace the observable $a$ in $\langle \Op_k(a) u_k,u_k\rangle_{L^2}$ by $a_k := a \chi_k$ where $\chi_k$ is a bump function equal to $1$ on $\{p=E_c\}$, supported on $\{p \in (E_c-k^{-\delta},E_c+k^{-\delta})\}$ for some $\delta < 1/2$. This allows to apply Egorov in times $t \simeq k^{\delta/2}$ as the exponential term in \eqref{equation:good} is then bounded by
\[
\exp\left((2(E-E_c)_+)^{1/2}tn\right) \leq \exp\left(\sqrt{2} k^{-\delta/2} k^{\delta/2} n\right) = \mc{O}(1).
\]
That is we obtain a similar estimate to \eqref{equation:parabolic}. Consequently, the following sequence of equalities hold for all $t \leq k^{\delta/2}$, modulo polynomial remainders (here $\delta > 0$ denotes a generic exponent which may differ from line to line):
\[
\begin{split}
\langle \Op_k(a) u_k,u_k\rangle_{L^2} & \overset{\eqref{equation:eigenstates}}{=} \langle e^{itk^{-1}\Delta_k}\Op_k(a)e^{-itk^{-1}\Delta_k} u_k,u_k\rangle_{L^2}  \\
& \overset{\eqref{equation:egorov}}{=} \langle \Op_k(a_k \circ \Phi_t) u_k,u_k\rangle_{L^2} + \mc{O}(k^{-\delta}) \\
& = \langle \Op_k\left(T^{-1}\int_0^T a_k \circ \Phi_t d t\right) u_k,u_k\rangle_{L^2} + \mc{O}(k^{-\delta}) \\
& \overset{\eqref{equation:burger}}{=} \int_{\{p=E_c\}} a(v) d\mu_{\mathrm{Liouv}}(v) \underbrace{\|u_k\|^2_{L^2}}_{=1}+ \mc{O}(k^{-\delta}),
\end{split}
\]
thus proving the claim.

\section{Hörmander's $L^\infty$ bound and magnetic zonal states}

\subsection{$L^\infty$ norm} A general result due to Hörmander \cite{Hormander-68} states that for
elliptic second-order differential operators on an $n$-dimensional compact manifold, one has the bound
\[
\|u_k\|_{L^\infty} \leq C k^{-(n-1)/2}\|u_k\|_{L^2},
\]
where $C>0$ is independent of $k \geq 0$ and $u_k$ is an eigenstate of the operator for an eigenvalue $\simeq k^2$. We also refer to the earlier works of Levitan and Avakumović
\cite{Levitan-52,Avakumovic-56} for related estimates. In the present
setting, the base manifold $\Sigma$ is two-dimensional, so the Hörmander estimate specializes to
\[
\|u_k\|_{L^\infty} \leq C k^{-1/2}\|u_k\|_{L^2}.
\]
In \cite[Theorem 1.1]{Chabert-Lefeuvre-26}, we establish the following dichotomy:

\begin{theorem}[$L^\infty$ bounds for magnetic eigenfunctions]
\label{theorem:Linfty}
The following assertions hold.
\begin{enumerate}[label=\emph{(\roman*)}]
\item \textbf{\emph{Low-energy regime.}}  
Assume that $0 \leq E < E_c$. Then there exists a sequence of eigenfunctions
$(u_k)_{k\geq0}$ solving \eqref{equation:eigenstates} such that
\begin{equation}
\label{equation:saturation}
\liminf_{k\to+\infty} k^{-1/2}
\|u_k\|_{L^\infty(\Sigma,L^{\otimes k})} > 0.
\end{equation}

\item \textbf{\emph{Critical energy regime.}}  
Assume that~\eqref{equation:eigenstates} holds with $E=E_c$ and that
$\varepsilon_k\leq h^\ell$ for some $0 < \ell \leq 1/15$. Then there exists a constant $C>0$ such that
\begin{equation}
\label{equation:polynomial-improvement}
\|u_k\|_{L^\infty(\Sigma,L^{\otimes k})}
\leq C k^{1/2 - \theta \ell/155800}.
\end{equation}

\end{enumerate}
\end{theorem}

In the high-energy regime ($E > E_c$), the magnetic flow is hyperbolic. It should then follow from now-standard arguments that, similarly to eigenfunctions of the Laplacian on functions, one has a Bérard-type $\sqrt{\log}$-improvement, that is:
\[
\|u_k\|_{L^\infty(\Sigma,L^{\otimes k})} \leq C (\log k)^{-1/2} k^{1/2}.
\]
See \cite{Berard-77,Bonthonneau-17} for further discussion.

In the low-energy regime, the eigenfunctions achieving the optimal growth rate predicted by
Hörmander exhibit a striking similarity with the zonal harmonics on the sphere. They may in fact be
constructed using closely related ideas, and their associated semiclassical defect measures can be
described explicitly; see \S\ref{ssection:zonal}. Motivated by this analogy, we refer to these
eigenfunctions as \emph{magnetic zonal states}.

The polynomial gain appearing in \eqref{equation:polynomial-improvement} at the critical energy
level is certainly not sharp, and no attempt was made to optimize the exponent. Determining the
optimal rate remains an open problem. It is actually a straightforward consequence of Theorem \ref{theorem:horocyclic} once combined with the following general bound for semiclassical operators established in \cite[Proposition 3.1]{Chabert-Lefeuvre-26}:

\begin{proposition}
\label{proposition:sup}
Let $P_h \in \Psi^{\mathrm{comp}}_h(\Sigma)$ be a semiclassical pseudodifferential operator and consider $(u_h)_{h > 0}$ such that
\[
P_h u_h = (E + o(1))u_h, \qquad \|u_h\|_{L^2}=1.
\]
Further assume that $p^{-1}(E-\delta,E+\delta)$ is compact for some $\delta > 0$, and $d_\xi p \neq 0$ on $p^{-1}(E-\delta,E+\delta)$. Then for all $0 \leq \eps < 1/2$, there exists $C > 0$ such that for all $x \in M$:
\[
|u_h(x)|^2 \leq  C \left(h^{-(n-1)-\eps} \int_{B(x,C h^\eps)} |u_h(x)|^2 + h^{-(n-1)+\eps}\right).
\]
\end{proposition}

This estimate was originally established in \cite{Donnelly-01} for the Riemannian Laplacian, in a slightly more precise form, following an argument due to Bourgain. To obtain the polynomial improvement (Theorem \ref{theorem:Linfty}, item (ii)), it suffices to combine Proposition \ref{proposition:sup} with Theorem \ref{theorem:horocyclic}, applied with the function $a = \pi^*\chi$, where $\pi : T^*\Sigma \to \Sigma$ is the projection and $\chi \in C^\infty(\Sigma)$ is a bump function localized on $B(x,Ch^\eps)$ (here $h=1/k$).

To the best of our knowledge, this phenomenon provides the first example
where the asymptotic behavior of $L^\infty$ norms for eigenfunctions of an elliptic operator
undergoes such a qualitative change depending on the energy range. In particular, it
would be interesting to understand the transition regime as the energy approaches $E_c$ from
below.

Sharper versions of Hörmander’s estimate have been studied in various geometric settings, notably
on manifolds of negative curvature; see for instance
\cite{Berard-77,Bonthonneau-17,Canzani-Galkowski-23}. Polynomial improvements, while conjectured in
several contexts such as negatively curved surfaces \cite{Sarnak-95}, are notoriously difficult to
establish and are currently known only in special cases, for example in
\cite{Iwaniec-Sarnak-95,Ingremeau-Vogel-24}.

\subsection{Magnetic zonal states}
\label{ssection:zonal}

We now describe the construction of the magnetic zonal states in the energy range
$0 < E < E_c$, which are responsible for the saturation of Hörmander’s estimate
in~\eqref{equation:saturation}.

Fix a point $x \in \Sigma$. The corresponding energy shell
in the cotangent fiber at $x$,
\[
C(x,E) := \{ p = E \} \cap T^*_{x}\Sigma,
\]
is a one-dimensional manifold diffeomorphic to a circle (in the degenerate case $E=0$,
the set $C(x,0)$ reduces to a single point and requires a separate argument).

Since the magnetic Hamiltonian flow $(\Phi_t)_{t\in\R}$ is periodic on the energy hypersurface
$\{p=E\}$, with period $T_E = (B^2 - 2E)^{-1/2}$, each point $(x,\xi)$ with $\xi \in C(x,E)$ generates a closed bicharacteristic:
\[
\Phi_{T_E}(x,\xi) = (x,\xi).
\]
Projecting this orbit onto the base manifold via the canonical projection
$\pi : T^*\Sigma \to \Sigma$ yields a smooth closed curve
$\gamma : [0,T_E] \to \Sigma$.

The magnetic zonal state is then defined by averaging these Gaussian beams over the
entire energy circle $C(x,E)$. More precisely, fixing an arbitrary
$\xi \in C(x,E)$, we set
\begin{equation}
\label{equation:magnetic-zonal-state}
\mathbf{u}_{k,x}
:= k^{1/4} \int_0^{2\pi}
\mathbf{f}_{k,x,R_\theta \xi}\,d\theta,
\end{equation}
where $R_\theta : T^*\Sigma \to T^*\Sigma$ is the rotation by angle $\theta$ in the
cotangent fibers. It can be verified that $1/C \leq \|\mathbf{u}_{k,x}\|_{L^2} \leq C$ for some uniform constant $C > 1$. The sequence $u_k := \mathbf{u}_{k,x}$ is precisely the family of eigenfunctions that
achieves the maximal $L^\infty$ growth described in
Theorem~\ref{theorem:Linfty}\, item (i).

However, it is an open question to determine
\[
\Omega := \limsup_{k\to+\infty} k^{-1/2}
\|u_k\|_{L^\infty(\Sigma,L^{\otimes k})} > 0,
\]
where $u_k \in C^\infty(\Sigma,L^{\otimes k})$ is such that $k^{-2}\Delta_k = (E+o(1)) u_k$ for $0 \leq E < E_c$ and $\|u_k\|_{L^2(\Sigma,L^{\otimes k})}=1$. On the sphere, the constant $\Omega$ is explicit, and equal to $(2\pi)^{-1/2}$. Finally, let us mention that there should exist an algebraic construction of such magnetic zonal states using the representation theory of $\mathrm{SL}(2,\R)$; this is left for future investigation.

\subsection{Defect measures of magnetic zonal states}

We now turn to the description of the semiclassical defect measure associated with the
sequence of magnetic zonal states~\eqref{equation:magnetic-zonal-state}. Microlocally, it can be easily seen that $(u_k)_{k\geq0}$ concentrates on the
two-dimensional invariant torus
\[
\begin{aligned}
\T^2(x,E)
&:= \bigl\{ \Phi_t(x,R_\theta\xi)
\;:\; \theta \in [0,2\pi],\; t \in [0,T_E] \bigr\} \\
&\simeq \R_t/(T_E\Z) \times \R_\theta/(2\pi\Z),
\end{aligned}
\]
where $\xi \in C(x,E)$, which is foliated by periodic magnetic trajectories. We denote by
$\Leb_{\T^2} := d\theta \otimes d t$ the Lebesgue measure on $\T^2(x,E)$. The following result is proved in \cite[Theorem 1.2]{Chabert-Lefeuvre-26}

\begin{theorem}[Semiclassical defect measure of magnetic zonal states]
Let $(u_k)_{k\geq 0}$ be the $L^2$-normalized sequence of magnetic zonal states defined in
\eqref{equation:magnetic-zonal-state}. Then:
\begin{enumerate}[label=\emph{(\roman*)}]
\item The semiclassical defect measure $\mu$
is supported on $\T^2(x,E)$ and equal to the normalized Lebesgue measure $\mu = (2\pi T_E)^{-1} d \theta \otimes d t$.
\item Its projection $\nu := \pi_* \mu$ is absolutely continuous with respect to the Riemannian volume $\vol_\Sigma$, that is $\nu = \alpha \vol_{\Sigma}$ for some explicit, upper semi-continuous function $\alpha \in L^1(\Sigma,\vol_{\Sigma})$. In addition, near $x$, one has
\[
\alpha(y) \sim_{y \to x} \dfrac{1}{4\pi T_E E d(x,y)}.
\]
\end{enumerate}
\end{theorem}

The measure $\nu$ has full support on $\Sigma$ provided $E \gg 0$ is large enough. This situation is closely analogous to that of zonal harmonics on the
sphere, whose semiclassical measures are supported on invariant tori associated with
closed geodesics. However, the density $\alpha$ only blows up at $y=x$ (whereas it blows up at two antipodal points on the sphere which correspond to conjugate points for the geodesic flow). This is due to the absence of ``magnetic'' conjugate points in this context.

\bibliographystyle{alpha}
\bibliography{biblio}

\end{document}